\documentclass[12pt,a4paper]{article}
\usepackage{amsmath,amsfonts,amssymb,amsthm,amscd}
\newtheorem{thm}{Theorem}
\newtheorem{lem}[thm]{Lemma}
\newtheorem{prop}[thm]{Proposition}

\newtheorem{cor}[thm]{Corollary}

\begin{document}

\title{Quaternionic connections, induced holomorphic structures
and a vanishing theorem}

\author{LIANA DAVID}

\maketitle

{\bf Abstract:} We classify the holomorphic structures of the
tangent vertical bundle $\Theta$ of the twistor fibration of a
quaternionic manifold $(M,Q)$ of dimension $4n\geq 8$. In
particular, we show that any self-dual quaternionic connection $D$
of $(M,Q)$ induces an holomorphic structure $\bar{\partial}^{D}$
on $\Theta .$ We construct a Penrose transform which identifies
solutions of the Penrose operator $P^{D}$ on $(M,Q)$ defined by $D$ with
the space of $\bar{\partial}^{D}$-holomorphic purely imaginary
sections of $\Theta .$ We 
prove that the tensor powers $\Theta^{s}$ (for any $s\in\mathbb{N}\setminus \{ 0\}$)
have no global non-trivial $\bar{\partial}^{D}$-holomorphic sections, when  
$(M, Q)$ is compact and has a compatible 
quaternionic-K\"{a}hler metric of negative (respectively,  
zero) scalar curvature
and the quaternionic connection $D$ is 
closed (respectively, closed but not exact).  
\footnote{{\bf Affiliation:} 
{\it Present:} Centro di Ricerca
Matematica ``Ennio de Giorgi'', Scuola Normale Superiore,
Piazza dei Cavalieri nr. 3, Pisa, Italia; e-mail:
l.david@sns.it; {\it Permanent:} Institute of Mathematics of the
Romanian Academy ``Simion Stoilow'', Calea Grivitei nr. 21,
Sector 1, Bucharest, Romania; e-mail: liana.david@imar.ro}\\

{\bf Key words and phrases:} quaternionic manifolds,
twistor spaces, Penrose transforms, Penrose operators.

{\bf Mathematical Subject Classification:} 53C26, 53C28, 53C15.\\

\section{Introduction}

An almost quaternionic structure on a manifold $M$ of dimension
$4n\geq 8$ is a rank three subbundle $Q$ of $\mathrm{End}(TM)$
locally generated by three almost complex structures which satisfy
the quaternionic relations. The bundle $Q$ has a natural Euclidian
metric, with respect to which any such system of almost complex
structures is orthonormal. $Q$ is called a quaternionic structure
if it is preserved by a torsion-free linear connection on $M$,
called a quaternionic connection. A quaternionic manifold is a
manifold together with a fixed quaternionic structure.

One of the main techniques to study quaternionic manifolds is
provided by twistor theory. The twistor space $Z$ of a
quaternionic manifold $(M,Q)$ is the total space of the unit
sphere bundle of $Q$, or the set of complex structures of tangent
spaces of $M$ which belong to $Q$. It has a natural integrable
almost complex structure which makes $Z$ a complex manifold. Of
interest in this paper is the tangent vertical bundle $\Theta$ of
the twistor fibration $\pi :Z\to M$, which is a hermitian
complex line bundle over $Z$. We show that any quaternionic
connection $D$ on $(M,Q)$ defines a hermitian connection $\nabla$
on $\Theta$, which is a Chern connection if and only if $D$ is
self-dual, i.e. the curvature of connection on
$\Lambda^{4n}(T^{*}M)$ induced by $D$ is $Q$-hermitian (see
Proposition \ref{l1}). When $D$ is  self-dual, the $(0,1)$-part
$\bar{\partial}^{D}$ of $\nabla$ is a real holomorphic structure
of $\Theta$, where by "real" we mean that the space of
$\bar{\partial}^{D}$-holomorphic sections of $\Theta$ is invariant
under the canonical anti-holomorphic involution of $Z$, defined as
the antipodal map along the fibers of $\pi$ (lifted to $\Theta$).
However, the complex line bundle $\Theta$ admits holomorphic
structures which are not necessarily real. Our first main result
is Theorem \ref{main} of Section \ref{proof}, and represents a
classification of all holomorphic structures of $\Theta$, in terms
of self-dual quaternionic connections on $(M,Q)$ and $1$-forms on
$M$ with $Q$-hermitian exterior derivative. This result is
analogous to Theorem $1$ of \cite{paul}, which classifies the
holomorphic structures of the tangent vertical bundle of the
twistor fibration of a conformal self-dual $4$-manifold, in terms
of self-dual Weyl connections and Maxwell fields on the conformal
$4$-manifold. As shown in \cite{paul}, the tangent vertical bundle
of the twistor fibration of a conformal self-dual $4$-manifold has
a canonical class of equivalent holomorphic structures, defined by
the Levi-Civita connections of the metrics in the conformal class.
Corollary \ref{c1} of Section \ref{proof} represents a similar
result in the quaternionic context. In the last two sections --
Section \ref{spenrose} and Section \ref{svanishing} -- we turn our
attention to the holomorphic sections of $\Theta$ (and, more generally,
of its positive tensor powers) with respect to
the holomorphic structures $\bar{\partial}^{D}$. More precisely,
in Section \ref{spenrose} we construct a Penrose transform, which
identifies the $\bar{\partial}^{D}$-holomorphic purely imaginary
sections of $\Theta$ with the kernel of the Penrose operator
$P^{D}$ of $(M,Q)$ defined by a self-dual quaternionic connection
$D$ (see Proposition \ref{lipsa}). 
In Section \ref{svanishing} we
prove that the Penrose operator $P^{D}$ has no global non-trivial
solutions when $(M,Q)$ is compact and admits a compatible
quaternionic-K\"{a}hler metric $g$ of negative or zero scalar curvature,
and $D$ is related to the Levi-Civita connection of $g$ in a
suitable way (see Theorem \ref{doi}). 
As a consequence of Theorem \ref{doi} and of the Penrose
transform, we deduce  that if $(M, Q, g)$ is compact 
quaternionic-K\"{a}hler 
and $D$ is a quaternionic connection 
on $(M, Q)$ such that $D$ is closed and $\mathrm{Scal}^{g}<0$
(respectively,  
$D$ is closed but not exact and
$\mathrm{Scal}^{g}=0$), then 
$\Theta$ has no global non-trivial $\bar{\partial}^{D}$-holomorphic sections  
(see Corollary \ref{co}). Similar results hold for all positive tensor powers 
of $\Theta$ (see the end of Section \ref{svanishing}). In the context
of conformal self-dual $4$-manifolds, this theory has been developed in
\cite{paul}.
Corollary \ref{co} has been proved in \cite{ponte}, using different 
methods. Other Penrose transforms and 
vanishing theorems on quaternionic-K\"{a}hler manifolds
have been developed in \cite{sal} and \cite{nitta}, where $\Theta$
(eventually coupled with the pull-back of a Yang-Mills field on $M$)
was considered with its canonical class of equivalent holomorphic
structures. In Sections \ref{spenrose} and \ref{svanishing} 
we generalise some results 
of \cite{nitta}, by considering $\Theta$ endowed with an 
holomorphic structure which does not necessarily 
belong to the canonical class of holomorphic structures of $\Theta .$

\section{Basic facts on quaternionic manifolds}\label{basic}

In this preliminary section we recall some basic facts we shall
need about quaternionic manifolds (quaternionic connections and
twistor spaces of quaternionic manifolds). We follow the treatment
of \cite{alek} and \cite{besse} (Chapter 14, Section $G$). All our
quaternionic manifolds will be connected and of dimension $4n\geq 8.$ For a
manifold $M$, $TM$, $T^{*}M$  and $\Omega^{k}(M)$ will denote the
real tangent bundle of $M$, the real cotangent bundle of $M$ and
the space of smooth real-valued $k$-forms on $M$, respectively.
For a vector bundle $V\to M$, $\Omega^{k}(M,V)$ will
denote the space of smooth $k$-forms on $M$ with values in $V$.
In our conventions, the curvature $R^{\nabla}$ of a connection
$\nabla$ acting on $V$ is defined by
$R^{\nabla}_{X, Y}s:=\nabla_{X}\nabla_{Y}s-\nabla_{Y}\nabla_{X}s-\nabla_{[X, Y]}s$, where 
$X, Y$ are vector fields on $M$ and $s$ is a smooth section of $V$.

\subsection{Quaternionic connections}\label{basicquat} Let $(M,Q)$ be a
quaternionic manifold and $D$ a quaternionic connection on
$(M,Q).$ Any other quaternionic connection $D^{\prime}$ is related
to $D$ by $D^{\prime}_{X}=D_{X}+S_{X}^{\alpha}$, where
$\alpha\in\Omega^{1}(M)$ and
\begin{equation}\label{a2}
S_{X}^{\alpha}:=\alpha (X)\mathrm{Id}_{TM}+\alpha\otimes X
-\sum_{i=1}^{3}[\alpha (J_{i}X)J_{i}+
(\alpha\circ J_{i})\otimes J_{i}X],\quad X\in TM.
\end{equation}
Here ``$\mathrm{Id}_{TM}$''denotes the identity endomorphism of
$TM$ and $\{ J_{1},J_{2},J_{3}\}$ is an admissible basis of $Q$,
i.e. a system of locally defined almost complex structures which
satisfy the quaternionic relations and generate $Q.$ The
connections $D$ and $D^{\prime}$ are {\it equivalent} if $\alpha$
is an exact $1$-form. We say that $D$ is {\it closed}
(respectively, {\it exact}) if it induces a flat connection on
$\Lambda^{4n}(T^{*}M)$ (respectively, if there is a volume form on
$M$ preserved by $D$). There always exist exact quaternionic
connections on $(M,Q)$: using relation  (\ref{a2}), one can check
that $D^{\prime}_{X}\underline{\mathrm{vol}}=
D_{X}\underline{\mathrm{vol}}-4(n+1)\alpha(X)\underline{\mathrm{vol}}$,
where $\underline{\mathrm{vol}}$ is an arbitrary volume form on
$M$; if $D_{X}\underline{\mathrm{vol}}=
\omega(X)\underline{\mathrm{vol}}$ for a $1$-form
$\omega\in\Omega^{1}(M)$, then the quaternionic connection
$D^{\prime}:=D+S^{\alpha}$ with $\alpha :=\frac{1}{4(n+1)}\omega$
is exact, because $\underline{\mathrm{vol}}$ is
$D^{\prime}$-parallel. Equally easy can be shown that any two
exact quaternionic connections are equivalent. The family of exact
quaternionic connections forms the {\it canonical class of
equivalent quaternionic connections of $(M,Q).$} We shall meet a
third class of connections, the so called {\it self-dual}
quaternionic connections; a quaternionic connection is self-dual
if the induced connection on the bundle $\Lambda^{4n}(T^{*}M)$ has
$Q$-hermitian curvature, i.e. its curvature is invariant 
with respect to any 
complex structure which belongs to $Q$.

A quaternionic curvature tensor of $(M,Q)$ is a curvature tensor
$R$ of $M$ (i.e. a section of
$\Lambda^{2}(T^{*}M)\otimes\mathrm{End}(TM)$ in the kernel of the
Bianchi map) which takes values in the normalizer of $Q$, i.e. for
any $X,Y\in TM$, $[R_{X,Y},Q]\subset Q.$ The space
$\mathcal{R}(N(Q))$ of quaternionic curvature tensors decomposes
into the direct sum $\mathcal{W}\oplus{\mathcal R}^{\mathrm{Bil}}$
where $\mathcal{W}$, called the space of quaternionic Weyl
curvatures, is the kernel of the Ricci contraction
$\mathrm{Ricci}: \mathcal{R}(N(Q))\to\mathrm{Bil}(TM)$,
defined by $\mathrm{Ricci}(R)_{X,Y}:=\mathrm{trace} \{
Z\to R_{Z,X}Y \} $ and $\mathcal{R}^{\mathrm{Bil}}$ is
isomorphic to the space $\mathrm{Bil}(TM)$ of bilinear forms on
$TM$, by means of  the isomorphism which associates to $\eta
\in\mathrm{Bil}(TM)$ the quaternionic curvature
\begin{align*}
R^{\eta }_{X,Y}:=&\left( \eta(Y,X)-\eta (X,Y)\right)\mathrm{Id}_{TM}-\eta_{X}\otimes Y
+\eta_{Y}\otimes X\\
&-\sum_{i=1}^{3}[(\eta(Y,J_{i}X)-\eta (X,J_{i}Y))J_{i}+
(\eta_{Y}\circ J_{i})\otimes J_{i}X
-(\eta_{X}\circ J_{i})\otimes J_{i}Y],
\end{align*}
where $X,Y\in TM$ and $\eta_{X}:=\eta (X,\cdot )$, $\eta_{Y}:=\eta (Y,\cdot ).$
Hence
the curvature of any quaternionic connection $D$ decomposes as
$R^{D}=W+R^{\eta}$, where $W\in \mathcal W$, called the quaternionic
Weyl tensor,
is an invariant of the quaternionic structure (i.e. is
independent of the choice of quaternionic connection) and
satisfies
$$
[W_{X,Y},A]=0,\quad X,Y\in TM,\quad A\in Q .
$$
With respect to an admissible basis $\{ J_{1},J_{2},J_{3}\}$ of
$Q$,
\begin{equation}\label{alta1}
[R^{D}_{X,Y},J_{i}]=\Omega_{k}(X,Y)J_{j}-\Omega_{j}(X,Y)J_{k},
\end{equation}
where $(i,j ,k)$ is a cyclic permutation of $(1,2,3)$ and, for any
$i \in \{ 1,2,3\}$ and $X, Y\in TM$,  
\begin{equation}\label{alta2}
\Omega_{i}(X, Y)=-\frac{1}{2n} \mathrm{tr}\left( J_{i}R^{D}_{X,
Y}\right) = 2\left( \eta (X,J_{i}Y)-\eta (Y,J_{i}X)\right) . 
\end{equation}

The bilinear form $\eta$ is related to the Ricci tensor $\mathrm{Ricci}(R^{D})$
of $D$ in the following way (see \cite{alek}, p. 223)
\begin{equation}\label{alta}
\eta =\frac{1}{4(n+1)}\mathrm{Ricci}(R^{D})^{\mathrm{skew}}+\frac{1}{4n}
\mathrm{Ricci}(R^{D})^{\mathrm{sym}}
-\frac{1}{2n(n+2)}P_{h}\left( \mathrm{Ricci}(R^{D})^{\mathrm{sym}}\right) ,
\end{equation}
where ``$\mathrm{sym}$'' and ``$\mathrm{skew}$''
denote the symmetric,
respectively skew-symmetric parts of a bilinear form and $P_{h}$
is the projection
$$
P_{h}(\eta ):=\frac{1}{4}\left( \eta +
\sum_{i=1}^{3}\eta (J_{i}\cdot ,J_{i}\cdot )\right)
$$
of $\mathrm{Bil}(TM)$ onto the space of $Q$-hermitian bilinear
forms. We mention that in general, $\mathrm{Ricci}(R^{D})$ is not
symmetric. More precisely, $\mathrm{Ricci}(R^{D})^{\mathrm{skew}}$
is half of the curvature of the connection induced by $D$ on the
canonical bundle $\Lambda^{4n}(T^{*}M)$ (see \cite{alek}, p. 222
and p. 224), so that $D$ has symmetric Ricci tensor (respectively,
the skew part of the Ricci tensor of $D$ is $Q$-hermitian) if and
only if $D$ is a closed (respectively, self-dual) quaternionic
connection. Finally, we remark that if $D^{\prime}=D+S^{\alpha}$,
then the Ricci tensors of $D$ and $D^{\prime}$ are related by the
following formulas (see \cite{alek}, p. 263)
\begin{align*}
\mathrm{Ricci}(R^{D^{\prime}})^{\mathrm{sym}}&=
\mathrm{Ricci}(R^{D})^{\mathrm{sym}}+4n\left(\alpha\otimes\alpha -
\sum_{i=1}^{3}(\alpha\circ J_{i})\otimes (\alpha\circ J_{i})
-(D \alpha )^{\mathrm{sym}}\right)\\
&+8P_{h}\left(\alpha\otimes\alpha -
\sum_{i=1}^{3}(\alpha\circ J_{i})\otimes (\alpha\circ J_{i})
-(D\alpha )^{\mathrm{sym}}\right)\\
\mathrm{Ricci}(R^{D^{\prime}})^{\mathrm{skew}}&=
\mathrm{Ricci}(R^{D})^{\mathrm{skew}}-4(n+1)
d\alpha.
\end{align*}

\subsection{Twistor theory of quaternionic manifolds}

\textbf{The twistor space of a quaternionic manifold:} As
mentioned in the Introduction, the twistor space $Z$ of $(M,Q)$,
defined as the total space of the unit sphere bundle of $Q$, 
has a natural complex
structure. In order to define it, we first consider a twistor line
$Z_{p}$, i.e the fiber of the natural projection $\pi
:Z\to M$ corresponding to a point $p\in M.$ Let
$\langle\cdot ,\cdot\rangle$ be the natural Euclidian metric of
the bundle $Q$. Then $T_{J}Z_{p}$ consists of all $J$-anti-linear
endomorphisms of $T_{p}M$ which belong $Q_{p}$, or to the
orthogonal complement $J^{\perp}$ of $J$ in $Q_{p}$, with respect
to the metric $\langle\cdot,\cdot\rangle$ . Note that $Z_{p}$ is a
K\"{a}hler manifold: it has a complex structure $\mathcal J$,
defined as ${\mathcal J}(A):=J\circ A$, for any $A\in T_{J}Z_{p}$,
and a compatible Riemannian metric, induced from the metric of
$Q_{p}$, since $T_{J}Z_{p}\subset Q_{p}$. Now we are able to
define the complex structure $\mathcal J$ of $Z$: chose a
quaternionic connection $D$ of $(M,Q)$. Since $D$ preserves $Q$ and
$\langle\cdot ,\cdot\rangle$, it induces a connection on the
twistor bundle $\pi :Z\to M$, i.e. a decomposition of
every tangent space $T_{J}Z$ into the vertical tangent space
$T_{J}Z_{p}$ and horisontal space $\mathrm{Hor}_{J}.$ On
$\mathrm{Hor}_{J}$, identified with $T_{p}M$ by means of the
differential $\pi_{*}$, $\mathcal{J}$ is equal to $J$. On
$T_{J}Z_{p}$, $\mathcal{J}$ is defined as above. It can be shown
that $\mathcal J$ so defined is independent of the choice of
quaternionic connection and is integrable. 
The twistor space $Z$ becomes a complex manifold of dimension $2n+1$ and
the twistor lines are
complex projective lines of $Z$ with normal bundle
$\mathbb{C}^{2n}\otimes{\mathcal O}(1)$.\\

\textbf{The tangent vertical bundle $\Theta$:} The tangent
vertical bundle $\Theta$ of the twistor projection $\pi
:Z\to M$ is the bundle over $Z$ whose fiber at a point 
$J\in Z$ is the tangent space $T_{J}Z_{p}$ at the twistor
line $Z_{p}$ defined by $\pi (J)=p.$ It is a  
complex line bundle over the complex
manifold $Z$, with complex structure of the fibers defined by the
complex structure of the twistor lines. Moreover, it has a
canonical hermitian metric $h(X,Y):=\frac{1}{2}\left( \langle
X,Y\rangle -i\langle{\mathcal J}X,Y\rangle\right)$, for any
$X,Y\in\Theta_{J}=T_{J}Z_{p}\subset Q_{p}.$ Due to this, there is
an isomorphism between Chern connections of $\Theta$ (i.e.
hermitian connections with $\mathcal J$-invariant curvature) and
holomorphic structures of $\Theta $, i.e. operators
$$
\bar{\partial}:\Gamma (\Theta )\to\Omega^{0,1}(Z,\Theta )
$$
which satisfy the Liebniz rule
$$
\bar{\partial}(fs)=f\bar{\partial}(s)+\bar{\partial}(f)s,\quad
f\in C^{\infty}(Z,\mathbb{C}),\quad s\in\Gamma (\Theta )
$$ and whose natural extension to the complex $\Omega^{0,*}(Z,\Theta )$ satisfies
$\bar{\partial}^{2}=0.$ The isomorphism associates to a Chern
connection $\nabla$ its $(0,1)$-part
$$
\bar{\partial}_{U}s:=\frac{1}{2}\left( \nabla_{U}s+{\mathcal
J}\nabla_{{\mathcal J}U}s\right) ,\quad U\in TZ,\quad s\in\Gamma
(\Theta ).
$$
Hence the study of holomorphic structures of $\Theta$ reduces to
the study of Chern connections.\\

\textbf{Distinguished sections of $\Theta$:} Note that any section
$A\in \Gamma (Q)$ defines a section $\tilde{A}$ of $\Theta$, by
the formula:
$$
\tilde{A}(J)=\Pi_{J} (A):=
\frac{1}{2}\left( A+J\circ A\circ J\right) =
A-\langle A,J\rangle J,\quad J\in Z,
$$
where the bundle homomorphism $\Pi: \pi^{*}Q\to \Theta$ is
the orthogonal projection onto $\Theta\subset \pi^{*}Q$ with
respect to the metric of $\pi^{*}Q$ induced by the natural
Euclidian metric $\langle\cdot ,\cdot \rangle$ of $Q$. Such
sections of $\Theta$ will be called {\it distinguished}. The
differential $\sigma_{*}:TZ\to TZ$ of the antipodal map
$\sigma :Z\to Z$, $\sigma (J)=-J$ induces an involution on
the space of smooth sections of $\Theta$, which associates to a
section $s$ the section $\bar{s}$ defined as follows: for any
$J\in Z$, $\bar{s}_{J}:= \sigma_{*}\left( s_{\sigma (J)}\right) .$
If $s:=\tilde{A}$ is distinguished, then $\bar{s}=-s$. This is why
the distinguished sections are also called {\it purely imaginary}.
Moreover, ${\mathcal J}s$ is {\it real}, i.e. $\overline{\mathcal
J s}={\mathcal J}s.$ The distinguished sections of $\Theta$ will
play a fundamental role in our treatment.

\section{Holomorphic structures on $\Theta$}\label{proof}

In this section, we adapt the arguments used in \cite{paul}
(Sections II2, II4, II5) to the quaternionic context. To keep our
text short, we refer to \cite{paul} whenever the analogy is
straightforward.

Consider a quaternionic connection $D$ on $(M,Q)$. Then $\pi^{*}D$
is a connection on the pull-back bundle $\pi^{*}Q$ and $\nabla :=
\Pi \circ \pi^{*}D$ is a connection on $\Theta .$ Since $D$
preserves $\langle\cdot ,\cdot \rangle$, the connection $\nabla$
preserves the Euclidian metric of $\Theta .$ Like in \cite{paul},
one shows that $\nabla$ is $\mathbb{C}$-linear, i.e. that $\nabla
{\mathcal J}=0$, where $\mathcal J$ denotes the complex structure
of the fibers of $\Theta .$

\begin{prop}\label{l1} The connection $\nabla$ is a Chern connection if
and only if $D$ is self-dual.
\end{prop}

\begin{proof}
For any distinguished section
$\tilde{A}$ of $\Theta$ and $U\in T_{J}Z$, with $\pi_{*}U=X\in T_{p}M$,
\begin{equation}\label{num}
\nabla_{U}\tilde{A}=\Pi_{J}\left(D_{X}A\right) - \langle J,A\rangle
v^{\bar{D}}(U).
\end{equation}
Here $\Pi_{J}\left( D_{X}A\right)$ 
is the orthogonal projection of $D_{X}A\in Q_{p}$
onto $J^{\perp}=\Theta_{J}$ and  
$v^{\bar{D}}(U)\in T_{J}Z_{p}=\Theta_{J}$ denotes the
vertical part of $U$ with respect to the connection $\bar{D}$
induced by $D$ on the twistor bundle $\pi :Z\to M.$ From
relation (\ref{num}) we obtain, like in \cite{paul} (see p. 585
and Appendix A), the following expression of the curvature
$R^{\nabla}:$
\begin{align*}
R^{\nabla}_{\tilde{X},\tilde{Y}}{A}&=\Pi_{J}\left( [R^{D}_{X,Y},A]\right)\\
R^{\nabla}_{B,C}{A}&=-\Omega_{p}(B,C){\mathcal J}(A)\\
R^{\nabla}_{\tilde{X},B}{A}&=0,
\end{align*}
where $\tilde{X} ,\tilde{Y}\in T_{J}Z$ are $\bar{D}$-horisontal
lifts of $X,Y\in T_{p}M$, $B,C\in T_{J}Z_{p}$, $A\in\Theta_{J}$,
$\Pi_{J}:Q_{p}\to J^{\perp}=\Theta_{J}$ is the orthogonal projection
and $\Omega_{p}$ is the K\"{a}hler form of the twistor line
$Z_{p}$, which is obviously $\mathcal J$-invariant. Hence $\nabla$
is a Chern connection if and only if the horisontal part of
$R^{\nabla}$ is $\mathcal J$-invariant, i.e. for every $J\in Z$
and $A\in Q$ with $A\perp J$,
\begin{equation}\label{conditie}
\Pi_{J}\left( [R^{D}_{JX\land JY-X\land Y},A]\right) =0.
\end{equation}

In order to study condition (\ref{conditie}), we take an
admissible basis $\{ J_{1},J_{2},J_{3}\}$ of $Q$ with $J=J_{1}$,
so that $A=\lambda_{2}J_{2}+\lambda_{3}J_{3}$ for some
$\lambda_{2} ,\lambda_{3}\in\mathbb{R}.$ Then $\Pi_{J}$ becomes
the projection onto the subspace generated by $J_{2}$ and $J_{3}.$
Recall now that $R^{D}=W+R^{\eta}$, for some
$\eta\in\mathrm{Bil}(TM)$ and that the quaternionic Weyl tensor
$W$ commutes with the endomorphisms of $Q$. Using relations
(\ref{alta1}) and (\ref{alta2}), we easily obtain:
\begin{align*}
\Pi_{J_{1}}\left( [R^{D}_{J_{1}X\land J_{1}Y-X\land Y}, A]\right)
&=\left( \Omega_{1}(J_{1}X,J_{1}Y)-\Omega_{1}(X,Y)\right) J_{1}A\\
&=-4\left(
\eta^{\mathrm{skew}}(J_{1}X,Y)+\eta^{\mathrm{skew}}(X,J_{1}Y)\right) J_{1}A.
\end{align*}
Using relation (\ref{alta}) we deduce that (\ref{conditie}) holds
if and only if 
$\mathrm{Ricci}(R^{D})^{\mathrm{skew}}$ is $Q$-hermitian,
i.e. $D$ is a self-dual quaternionic connection.

\end{proof}

Self-dual quaternionic connections exist on any quaternionic
manifold.  If $D$ is a self-dual quaternionic connection, then any
other self-dual quaternionic connection is of the form
$D^{\prime}=D+S^{\alpha}$, with $d\alpha$ $Q$-hermitian. The way
the Chern connections of $\Theta$ determined by two self-dual
quaternionic connections of $(M,Q)$ are related is described in
the following proposition.

\begin{prop}\label{l2} Let $D$, $D^{\prime}=D+S^{\alpha}$ be two
self-dual quaternionic connections. Then the Chern connections
$\nabla$ and $\nabla^{\prime}$ induced by $D$ and $D^{\prime}$ are
related as follows:
$$
\nabla^{\prime}=\nabla +2{\mathcal J}\left( \pi^{*}\alpha\right) \otimes \mathcal{J}.
$$
\end{prop}

\begin{proof}
Fix an arbitrary $J\in Z$ and an admissible basis $\{
J_{1},J_{2},J_{3}\}$ of $Q$ with $J=J_{1}.$ Any $A\in\Theta_{J}$
is of the form $\lambda_{2}J_{2}+\lambda_{3}J_{3}$, for some
$\lambda_{2},\lambda_{3}\in\mathbb{R}.$ For every $U\in T_{J}Z$
with $\pi_{*}U=X$, we see, from relation (\ref{num}), that
\begin{align*}
(\nabla^{\prime}-\nabla)_{U}(A)=\Pi_{J_{1}}[D^{\prime}_{X}-D_{X}, A]&
=2\alpha (J_{1}X)(-\lambda_{2}J_{3}+\lambda_{3}J_{2})\\
&=-2\alpha (J_{1}X)J_{1}(\lambda_{2}J_{2}+\lambda_{3}J_{3})\\
&=2({\mathcal J}\pi^{*}\alpha )(U)J_{1}A .
\end{align*}
\end{proof}

\textbf{Remark:} Recall that two holomorphic structures on a
complex line bundle $V\to N$ over a complex manifold
$(N,J)$ are equivalent, if they are conjugated by an element of
the gauge group $C^{\infty}(N,\mathbb{C}^{*})$. Suppose now that
$V$ has a hermitian structure. Let $\bar{\partial}^{1}$ and
$\bar{\partial}^{2}$ be two holomorphic structures of $V$ and let
$\nabla^{1}$ and $\nabla^{2}$ be the corresponding Chern connections, with
$(0,1)$-parts $\bar{\partial}^{1}$ and $\bar{\partial}^{2}$
respectively. Then $\bar{\partial}^{1}$ and $\bar{\partial}^{2}$
are equivalent if and only if

\begin{equation}\label{chern}
\nabla^{2}=\nabla^{1}+(d^{J}\mathrm{log}\rho -d\theta )\otimes{\mathcal J},
\end{equation}
where $\rho$ is a positive smooth function, $\theta$ is a smooth
function with values in $S^{1}$ and $\mathcal J$ is the complex
structure of the fibers of $V$. The connections $\nabla^{1}$ and
$\nabla^{2}$ are equivalent as hermitian connections if
$d^{J}\mathrm{log}\rho =0.$

The following classification theorem holds:

\begin{thm}\label{main} Any holomorphic structure of $\Theta$ is
equivalent with an holomorphic structure
$\bar{\partial}^{D,\beta}:=\bar{\partial}^{D}+\tilde{\beta}$,
where $\bar{\partial}^{D}$ is the $(0,1)$-part of the Chern
connection of $\Theta$ induced by a self-dual quaternionic
connection $D$ of $(M,Q)$, $\beta\in\Omega^{1}(M)$ has
$Q$-hermitian exterior differential and $\tilde{\beta}\in
\Omega^{0,1}(Z,\mathrm{End}_{\mathbb{C}}(\Theta ))$ is defined as
follows: for any $U\in T_{J}Z$ with $\pi_{*}U=X$ and $s\in\Gamma
(\Theta )$,
$$
\tilde{\beta}_{U}(s):=\frac{1}{2}\left(\beta (X){\mathcal
J}s-\beta (JX)s\right) .
$$
Moreover, two holomorphic structures $\bar{\partial}^{D,\beta}$
and $\bar{\partial}^{D^{\prime},\beta^{\prime}}$ are equivalent if
and only if $D$ and $D^{\prime}$ are equivalent as quaternionic
connections of $(M,Q)$ and $d^{\beta}:=d+i\beta$ and
$d^{\beta^{\prime}}:=d+i\beta^{\prime}$ are equivalent as
hermitian connections of the hermitian trivial line bundle
$M\times\mathbb{C}$.
\end{thm}

\begin{proof}
The proof follows the same steps as the proof of Theorem $1$ of
\cite{paul} (note that, our Proposition \ref{l1} corresponds to
Proposition $2$ of \cite{paul} and our Proposition \ref{l2}
corresponds to Lemma $4$ of \cite{paul}). Due to this, we content
ourselves to explain why $\bar{\partial}^{D,\beta}$ is an
holomorphic structure. Since $d\beta$ is $Q$-hermitian, the
pull-back connection $d+i\pi^{*}\beta$ is a Chern connection on
the hermitian trivial line bundle $Z\times\mathbb{C}$. Let
$\nabla$ be the Chern connection of $\Theta$ induced by a
self-dual quaternionic connection $D$, as in Proposition \ref{l1}.
The tensor product connection $\nabla^{\beta}:=\nabla\otimes
(d+i\pi^{*}\beta )= \nabla +\pi^{*}\beta \otimes{\mathcal J}$ on
$\Theta=\Theta\otimes_{\mathbb{C}}\mathbb{C}$ is also a Chern
connection on $\Theta$. It can be checked that its
$(0,1)$-part is precisely $\bar{\partial}^{D,\beta}$. In
particular, $\bar{\partial}^{D,\beta}$ is an holomorphic structure
of $\Theta .$

\end{proof}

Recall now that any two exact quaternionic connections are
equivalent. The following Corollary is a  consequence of Theorem
\ref{main}.

\begin{cor}\label{c1} The tangent vertical bundle of the twistor
fibration of a quaternionic manifold $(M,Q)$ has a canonical class
of equivalent holomorphic structures, determined by the exact
quaternionic connections of $(M,Q)$.

\end{cor}

\section{A Penrose transform}\label{spenrose}

We shall use the $E-H$ formalism developed in
\cite{salamon}. We begin with a brief review of some basic facts
we shall need about the representation theory of the group
$Sp(1).$ Let $H\cong\mathbb{C}^{2}$ be an abstract $2$-dimensional
complex vector space on which $Sp(1)\cong SU(2)$ acts by complex
linear transformations, leaving
invariant a complex symplectic form $\omega$ and a compatible
quaternionic structure, i.e. a $\mathbb{C}$-anti-linear map
$q:H\to H$, which satisfies $q^{2}=-\mathrm{Id}_{H}$,
$\omega (qv,qw)=\overline{{\omega}(v,w)}$ and $\omega (v,qv)>0$,
for any $v,w\in H.$ The $2$-form $\omega$ together with $q$ define
an invariant hermitian positive definite metric $\langle \cdot
,\cdot\rangle :=\omega (\cdot ,q\cdot )$ on $H$. By means of the
identification $H\ni h\to \omega (h,\cdot )\in H $ between
$H$ and its dual $H^{*}$, $S^{2}(H)\subset H\otimes H\cong
H^{*}\otimes H\subset\mathrm{End}(H)$ acts on $H$ and its real
part (with respect to the real structure induced by $q$) is
isomorphic to the Lie algebra $sp(1)\subset \mathrm{End}(H)$ of
imaginary quaternions. We also need to recall that $H\otimes
S^{2}(H)$ has two $Sp(1)$-irreducible components: $S^{3}(H)$,
which is the kernel of the map $F:{H}\otimes S^{2}({H})\to
{H}$ defined by
\begin{equation}\label{F}
F(h,h_{1}h_{2}+h_{2}h_{1})=\omega (h,h_{1})h_{2}+
\omega (h,h_{2})h_{1} ,\quad h_{1},h_{2},h\in {H},
\end{equation}
and $H$, isomorphic to the hermitian orthogonal of $S^{3}({H})$ in
${H}\otimes S^{2}({H})$ with respect to the hermitian metric of
${H}\otimes S^{2}({H})$ induced by the hermitian metric $\langle
\cdot ,\cdot \rangle$ of $H$ (to simplify notations, we 
omit the tensor product signs, so that $h_{1}h_{2}+h_{2}h_{1}$
denotes $h_{1}\otimes h_{2}+h_{2}\otimes h_{1}$).

Coming back to geometry, the quaternionic structure of $(M,Q)$
determines a $G=GL(n,\mathbb{H})Sp(1)$ structure $F_{0}$, i.e. a
$G$-subbundle of the principal frame bundle of $M$, consisting of
all frames $f:T_{p}M\to \mathbb{H}^{n}$ which convert the
standard basis of imaginary quaternions, acting by multiplication
on $\mathbb{H}^{n}$ on the right, onto an admissible basis of
$Q_{p}$, acting naturally on $T_{p}M$. Any representation of
$\tilde{G}=GL(n,\mathbb{H})\times Sp(1)$ determines a locally
defined bundle over $M$, which is globally defined when the
representation descends to $G$ (in which case the bundle is
associated to the principal $G$-bundle $F_{0}$). Real
representations and real vector bundles over $M$ will be
automatically complexified. There are two locally defined complex
vector bundles $\bf E$ and $\bf H$ on $M$, which are associated to
the standard representations of $GL(n,\mathbb{H})$ and $Sp(1)$ on
$E\cong\mathbb{C}^{2n}$ and $H\cong \mathbb{C}^{2}$ respectively,
extended trivially to $GL(n,\mathbb{H})\times Sp(1)$. The
$Sp(1)$-invariant structures of $H$ induce similar structures on
the bundle ${\bf H}$, which will be denoted with the same symbols
(e.g. $\omega$ will denote the symplectic form of $H$ as well as
the induced symplectic form on the bundle ${\bf H}$; in
particular, we shall identify ${\bf H}$ with its dual ${\bf
H}^{*}$ by means of the isomorphism ${\bf H}\ni h\to
\omega (h,\cdot)\in{\bf H}^{*}$; similarly, $\langle\cdot
,\cdot\rangle$ will denote the hermitian inner product of $H$ and
the induced hermitian metric on the bundle ${\bf H}$). Some of the
natural bundles over $M$ are isomorphic with tensor products and
direct sums of ${\bf H}$ and ${\bf E}$. For example, $TM$ is
isomorphic with ${\bf E}\otimes {\bf H}$, $T^{*}M$ with ${\bf
E}^{*}\otimes{\bf H}$, $Q$ with $S^{2}(\bf H)$ and the product
$T^{*}M\otimes Q$ decomposes as
\begin{equation}\label{d}
T^{*}M\otimes Q \cong {\bf E}^{*}\otimes{\bf H}\otimes S^{2}({\bf H})
\cong {\bf E}^{*}\otimes S^{3}({\bf H})\oplus {\bf E}^{*}
\otimes {\bf H},
\end{equation}
since $H\otimes S^{2}(H)\cong S^{3}(H)\oplus H$. The {\it Penrose
operator} $P^{D}:\Gamma (Q)\to \Gamma ({\bf E}^{*}\otimes
S^{3}({\bf H}))$ defined by a quaternionic connection $D$ of
$(M,Q)$ is the composition of $D :\Gamma (Q)\to \Gamma
(T^{*}M\otimes Q)$ with the projection onto the first component of
the decomposition (\ref{d}).

\begin{prop}\label{lipsa} Let $D$ be a self-dual quaternionic connection on
$(M,Q)$ and $A\in \Gamma (Q).$ Then the distinguished section
$\tilde{A}$ of $\Theta$ is $\bar{\partial}^{D}$-holomorphic if and
only if $A$ is a solution of the Penrose operator $P^{D}.$
\end{prop}

\begin{proof}
The section
$\tilde{A}$ is $\bar{\partial}^{D}$-holomorphic if and only if
it satisfies
\begin{equation}\label{ho}
\nabla_{\mathcal{J}U}(\tilde{A})={\mathcal J}\nabla_{U}(\tilde{A}),
\quad \forall U\in TZ,
\end{equation}
where $\nabla$ is the Chern connection of $\Theta$ induced
by $D$. Using relation (\ref{num}), it can be seen that (\ref{ho}) is equivalent
with
$$
D_{JX}A -\langle D_{JX}A,J\rangle J-JD_{X}A-\langle
D_{X}A,J\rangle\mathrm{Id}_{TM}=0,\quad \forall X\in
TM,\quad\forall J\in Z.
$$
For every unit $j\in sp(1)\subset \mathrm{End}(E\otimes H)$
(acting trivially on $E$), define
$$
T_{j}:E^{*}\otimes H^{*}\otimes S^{2}(H)\to E^{*}\otimes
H^{*}\otimes S^{2}(H)
$$
in the following way:  for any $\gamma\in E^{*}\otimes
H^{*}\otimes S^{2}(H)$ and $v\in E\otimes H$,
\begin{equation}\label{tj}
T_{j}(\gamma )(v):=\gamma (jv)-\langle \gamma (jv),j\rangle j -
j\circ \gamma(v) -\langle \gamma (v),j\rangle
\mathrm{Id}_{E\otimes H}.
\end{equation}
Here $\langle\cdot ,\cdot\rangle$  denotes the hermitian inner
product of $H\otimes H$ induced by the $Sp(1)$-invariant hermitian
inner product $\langle\cdot ,\cdot \rangle$ of $H$, i.e.
$$
\langle h_{1}h_{2},h_{3}h_{4}\rangle =\frac{1}{2}\langle
h_{1},h_{3}\rangle \langle h_{2},h_{4}\rangle ,\quad \forall
h_{i}\in H,
$$
so that the restriction of $\langle\cdot,\cdot\rangle$ to
$sp(1)\subset S^{2}(H)$ induces the natural Euclidian metric on
the bundle $Q$ (recall that $Q$ is associated to the adjoint
representation of $Sp(1)$ on its Lie algebra $sp(1)$, extended
trivially to $\tilde{G}$). The group $Sp(1)$ acts naturally on
$H^{*}\otimes S^{2}(H)\subset H^{*}\otimes\mathrm{End}(H)$, by
$$
(a \cdot\alpha )(h)=a\circ\alpha (a^{-1}h)\circ a^{-1}, \quad a\in
Sp(1),\quad \alpha\in H^{*}\otimes S^{2}(H),\quad h\in H.
$$
We can extend this action to an action of $\tilde{G}$ on
$E^{*}\otimes H^{*}\otimes S^{2}(H)$, with $GL(n,\mathbb{H})$
acting naturally on $E^{*}$. This extended action preserves the
$\mathbb{C}$-linear condition
\begin{equation}\label{lincon}
T_{j}(\gamma )=0,\quad \forall j\in sp(1), \quad
j^{2}=-\mathrm{Id}_{H}
\end{equation}
on $E^{*}\otimes H^{*}\otimes S^{2}(H)$, since
$$
T_{j}(a\cdot\gamma )(v)=a\circ T_{j^{\prime}}(\gamma
)(a^{-1}v)\circ a^{-1},\quad\forall a\in\tilde{G},\quad \forall v\in E\otimes H,
$$
where $j^{\prime}:=a^{-1}\circ j\circ a.$ Since $E^{*}\otimes H$
and $E^{*}\otimes S^{3}(H)$ are $\tilde{G}$-irreducible components of
$E^{*}\otimes H\otimes S^{2}(H)$ and there are distinguished
sections of $\Theta$ which are not
$\bar{\partial}^{D}$-holomorphic, from Shur's lemma it is enough
to check that any element of $E^{*}\otimes H$ satisfies
(\ref{lincon}). This can be done in the following way: let
$e^{*}h\in E^{*}\otimes H$ be decomposable. Without loss of
generality, we can take $\langle h,h\rangle =1$. Define
$\tilde{h}:=q(h).$ Then the basis $\{ h,\tilde{h}\}$ is unitary
with respect $\langle\cdot ,\cdot \rangle$ and $\omega
(h,\tilde{h})=1.$ As an element of $E^{*}\otimes H\otimes
S^{2}(H)\subset E^{*}\otimes H\otimes\mathrm{End}(E\otimes H)$,
$e^{*}h$ has the following form
\begin{equation}\label{alpha}
\gamma (v)=\left( 2 (e^{*}\tilde{h})(v)hh-
(e^{*}h)(v)(\tilde{h}h+h\tilde{h})\right)\mathrm{Id}_{E},
\quad\forall  v\in E\otimes H.
\end{equation}
The identity endomorphism of $H$
is $\mathrm{Id}_{H}=h\tilde{h}-\tilde{h}h$
and a basis of unit imaginary quaternions can be chosen to be
\begin{align*}
j_{1}&:=-i(h\tilde{h}+\tilde{h}h)\mathrm{Id}_{E}\\
j_{2}&:=-(hh+\tilde{h}\tilde{h})\mathrm{Id}_{E}\\
j_{3}&:=i(\tilde{h}\tilde{h}-hh)\mathrm{Id}_{E}.
\end{align*}

Consider now $j=a_{1}j_{1}+a_{2}j_{2}+a_{3}j_{3}\in sp(1)$
an arbitrary unit imaginary quaternion (so that
the $a_{i}$'s are real and
$a_{1}^{2}+a_{2}^{2}+a_{3}^{2}=1$).
Then, for every $v\in E\otimes H$,
\begin{align*}
\gamma (jv)&=2hh\left( a_{1}i(e^{*}\tilde{h})(v)+(a_{2}+ia_{3})(e^{*}h)(v)\right)
\mathrm{Id}_{E}\\
&+(h\tilde{h}+\tilde{h}h)\left( a_{1}i(e^{*}h)(v)+(a_{2}-
ia_{3})(e^{*}\tilde{h})(v)\right)\mathrm{Id}_{E};\\
\langle \gamma (jv),j\rangle &=(-a_{2}+ia_{3})\left( a_{1}i(e^{*}\tilde{h})(v)
+(a_{2}+ia_{3})(e^{*}h)(v)\right)\\
&+a_{1}\left( -a_{1}(e^{*}h)(v)+(a_{3}+ia_{2})(e^{*}\tilde{h})(v)\right);\\
j\circ\gamma (v)&=(e^{*}h)(v)\left( (a_{2}+ia_{3})hh +(-a_{2}+ia_{3})\tilde{h}\tilde{h}
+ia_{1}(h\tilde{h}-\tilde{h}h)\right)\mathrm{Id}_{E}\\
&+2(e^{*}\tilde{h})(v)\left( ia_{1}hh+(a_{2}-ia_{3})h\tilde{h}\right)
\mathrm{Id}_{E};\\
\langle \gamma (v),j\rangle &=(e^{*}\tilde{h})(v)(-a_{2}+ia_{3})-(e^{*}h)(v)ia_{1}.
\end{align*}
From these relations it is straightforward to check, using the definition
 of $T_{j}$ given in (\ref{tj}), that
$T_{j}(\gamma )=0.$

\end{proof}

\section{A vanishing theorem}\label{svanishing}

In this section we consider a quaternionic manifold $(M,Q)$ which
admits a compatible quaternionic-K\"{a}hler metric $g$, i.e. the
Levi-Civita connection $D^{g}$ of $g$ is a quaternionic connection
of $(M,Q)$ and the endomorphisms of $Q$ are skew-symmetric with
respect to $g$. Let $g^{*}:T^{*}M\to TM$ be the
isomorphism defined by $g$.
Borrowing the terminology of \cite{uwe}, we define
a conformal weight operator
$$
B:T^{*}M\otimes Q\to T^{*}M\otimes Q
$$
by the following formula:
\begin{equation}\label{B}
B(\alpha\otimes A)(X):=[S^{\alpha}_{X},A],\quad \forall \alpha\in T^{*}M,\quad\forall A\in Q,\quad
\forall X\in TM,
\end{equation}
where $S^{\alpha}$ was defined in (\ref{a2}).
(A conformal weight operator has been defined
in \cite{paul}, for vector bundles on conformal manifolds,
associated to the principal bundle of conformal frames, and in
\cite{uwe} for vector bundles associated to the reduced frame
bundle of a Riemannian manifold with special holonomy).
The following lemma is a straightforward calculation:

\begin{lem}\label{confop} Let $\{ J_{1},J_{2},J_{3}\}$ be an admissible
basis of $Q$. Then for every $\alpha \in T^{*}M$,
$A\in Q$ and $X\in TM$, 
$$
B(\alpha\otimes A)(X)=\alpha \left( [J_{1},A](X)\right) J_{1}
+\alpha \left( [J_{2},A](X)\right) J_{2}
+\alpha \left( [J_{3},A](X)\right) J_{3}.
$$
\end{lem}

\begin{prop}\label{imp} Consider the decomposition (\ref{d}) of
$T^{*}M\otimes Q.$ The conformal weight operator $B$ acts as
$-2\cdot\mathrm{Id}_{{\bf E}^{*} \otimes S^{3}({\bf H})}$ on ${\bf
E}^{*}\otimes S^{3}({\bf H})$ and as $4\cdot\mathrm{Id}_{{\bf
E}^{*}\otimes {\bf H}}$ on ${\bf E}^{*}\otimes {\bf H}.$
\end{prop}

\begin{proof}

Let $e^{*}h\in\Gamma ({\bf E}^{*}\otimes {\bf H})$ be a
decomposable local section, with $\langle h,h\rangle =1.$ Define
$\tilde{h}:=q(h).$ As a section of $T^{*}M\otimes Q$,
$e^{*}h=\alpha\otimes A+\tilde{\alpha}\otimes\tilde{A}$, where
\begin{align*}
A &:=2hh \mathrm{Id}_{\bf E}\in\Gamma (Q),\\
\tilde{A}&:=- (\tilde{h}h+h\tilde{h})\mathrm{Id}_{\bf E} \in\Gamma
( Q)\\
\alpha &:=e^{*}\tilde{h}\in\Omega^{1}(M)\\
\tilde{\alpha}&:=e^{*}h\in \Omega^{1}(M).
\end{align*}
As in the proof of Proposition \ref{lipsa}, we consider the
following admissible basis of $Q$:

\begin{align*}
J_{1}&:=-i(h\tilde{h}+\tilde{h}h)\mathrm{Id}_{\bf E}\\
J_{2}&:=-(hh+\tilde{h}\tilde{h})\mathrm{Id}_{\bf E}\\
J_{3}&:=i(\tilde{h}\tilde{h}-hh)\mathrm{Id}_{\bf E}.
\end{align*}
It is easy to check the equalities:
\begin{align*}
[J_{1},A]&=-2i[h\tilde{h}+\tilde{h}h,hh]\mathrm{Id}_{\bf E}
=4ihh\mathrm{Id}_{\bf E}\\
[J_{2},A]&=-2[hh+\tilde{h}\tilde{h},hh]\mathrm{Id}_{\bf E}
=2(h\tilde{h}+\tilde{h}h)\mathrm{Id}_{\bf E}\\
[J_{3},A]&=2i[\tilde{h}\tilde{h}-hh,hh]\mathrm{Id}_{\bf E}=
-2i(h\tilde{h}+\tilde{h}h)\mathrm{Id}_{\bf E}.
\end{align*}
Using Lemma \ref{confop}, we get that
$$
B(\alpha\otimes A)(X)=-4
(e^{*}h)(X)(h\tilde{h}+\tilde{h}h)\mathrm{Id}_{\bf E}
+4(e^{*}\tilde{h})(X)hh\mathrm{Id}_{\bf E},\quad X\in TM.
$$
A similar calculation shows that
$B(\tilde{\alpha}\otimes\tilde{A})(X)=4(e^{*}\tilde{h})(X)hh
\mathrm{Id}_{\bf E},$
which readily implies that
$$
B(\alpha\otimes A+\tilde{\alpha}\otimes\tilde{A})=4(\alpha\otimes A+\tilde{\alpha}\otimes\tilde{A}),
$$
i.e. that
${\bf E}^{*}\otimes{\bf H}$ is the eigenspace of
$B$ corresponding to the eigenvalue $4$.
In order to show that ${\bf E}^{*}\otimes S^{3}({\bf H})$ is
the eigenspace of $B$ corresponding to the eigenvalue $2$,
we notice that ${\bf E}^{*}\otimes S^{3}({\bf H})\subset T^{*}M\otimes Q$ is locally
generated by sections $\gamma_{i}^{e^{*}}$ (for
$e^{*}\in\Gamma ({\bf E}^{*})$ and $i\in\{ 1,\cdots 4\}$)
defined as follows: for every $X\in TM$,
\begin{align*}
 \gamma^{e^{*}}_{1}(X)&:=\left( (e^{*}h)(X)(h\tilde{h}+\tilde{h}h)
+(e^{*}\tilde{h})(X)hh\right)\mathrm{Id}_{\bf E}\\
\gamma^{e^{*}}_{2}(X)&:= \left( (e^{*}\tilde{h})(X)(h\tilde{h}+\tilde{h}h)
+(e^{*}h)(X)\tilde{h}\tilde{h}\right)\mathrm{Id}_{\bf E}\\
\gamma^{e^{*}}_{3}(X)&:=(e^{*}h)(X)hh\mathrm{Id}_{\bf E}\\
\gamma^{e^{*}}_{4}(X)&:=(e^{*}\tilde{h})(X)\tilde{h}\tilde{h}\mathrm{Id}_{\bf E}.
\end{align*}
As before, one checks that
$B(\gamma^{e^{*}}_{i})=-2\gamma^{e^{*}}_{i}$, for every
$e^{*}\in\Gamma ({\bf E}^{*})$ and
$i\in\{ 1,\cdots ,4\}$. The conclusion follows.

\end{proof}

\begin{prop}\label{propo} Let $D=D^{g}+S^{\alpha}$ be a quaternionic
connection on a quaternionic-K\"{a}hler manifold $(M,Q,g)$, with
$\alpha\in\Omega^{1}(M)$ a co-closed $1$-form. Let $A\in\Gamma
(Q)$ be a solution of the Penrose operator $P^{D}$. Then
$$
\langle\mathrm{trace}_{g}(D^{2}A),A\rangle =-2|A|^{2}
\left(\frac{1}{4(n+2)}\mathrm{Scal}^{g}-2|\alpha |^{2}\right) .
$$
\end{prop}

\begin{proof}
It is straightforward to check that $D\circ B=\tilde{B}\circ D$,
where the connection $D$ (on both sides of this equality) 
acts on
$T^{*}M\otimes Q$ and $\tilde{B} :=\mathrm{Id}_{T^{*}M}\otimes B$
is an automorphism of $T^{*}M\otimes T^{*}M\otimes Q.$ Define
$$
\mathrm{trace}_{g}(\tilde{B}):T^{*}M\otimes T^{*}M\otimes
Q\to Q
$$
as follows: for any $A\in Q$, $\alpha ,\beta\in T^{*}M$,
$$
\mathrm{trace}_{g}(\tilde{B})(\alpha\otimes\beta\otimes
A):=\sum_{i=1}^{4n} \tilde{B}(\alpha\otimes\beta\otimes
A)(e_{i},e_{i})=B(\beta\otimes A)(g^{*}\alpha ) ,
$$
where $\{ e_{1},\cdots ,e_{4n}\}$ is an arbitrary $g$-orthonormal
basis of $TM$.  
Writing $A=a_{1}J_{1}+a_{2}J_{2}+a_{3}J_{3}$ in
terms of an admissible basis $\{ J_{1},J_{2},J_{3}\}$ of $Q$, we
readily obtain, from Lemma \ref{confop}, that
\begin{align*}
\mathrm{trace}_{g}(\tilde{B})(D^{2}A)&=
\sum_{i<j}\left( g\left([J_{1},R^{D}_{e_{i},e_{j}}(A)]e_{i},e_{j}\right) J_{1}
+ g\left( [J_{2},R^{D}_{e_{i},e_{j}}(A)]e_{i},e_{j}\right) J_{2}\right)\\
&+\sum_{i<j} g\left( [J_{3},R^{D}_{e_{i},e_{j}}(A)]e_{i},e_{j}\right) J_{3},
\end{align*}
where $R^{D}_{e_{i},e_{j}}(A)=[R^{D}_{e_{i},e_{j}},A]$ is the
commutator of the endomorphisms $R^{D}_{e_{i},e_{j}}$ and $A$ of
$TM$. In particular,
\begin{align*}
\langle\mathrm{trace}_{g}(\tilde{B})(D^{2}A),A\rangle &=
\sum_{i<j}g\left( [A,R^{D}_{e_{i},e_{j}}(A)]e_{i},e_{j}\right)\\
&=\sum_{i<j}\sum_{k,p}a_{k}a_{p}g\left(
[J_{k},R^{D}_{e_{i},e_{j}}(J_{p})] e_{i},e_{j}\right) .
\end{align*}
Using relations (\ref{alta1}) and (\ref{alta2}) we readily get that,
for every $k\in\{ 1,2,3\}$,
$$
\sum_{i<j}g\left([J_{k},R^{D}_{e_{i},e_{j}}(J_{k})]e_{i},e_{j}\right)
=-8\mathrm{trace}_{g}( \eta )
$$
and for every $k\neq p$,
$$
\sum_{i<j}g\left( [J_{k},R^{D}_{e_{i},e_{j}}(J_{p})]e_{i} +[J_{p},
R^{D}_{e_{i},e_{j}}(J_{k})]e_{i},e_{j}\right) =0
$$
from where we conclude that
\begin{equation}\label{tilde}
\langle \mathrm{trace}_{g}(\tilde{B})(D^{2}A),A\rangle
=-8|A|^{2}\mathrm{trace}_{g}(\eta ).
\end{equation}
Recall now that $\eta$ is related to the Ricci curvature
$\mathrm{Ricci}(R^{D})$ of $D$ as in (\ref{alta}).
According to Section \ref{basic}, we can express
$\mathrm{Ricci}(R^{D})$ in terms of $\alpha$ and the Ricci tensor
of $g$, so that, taking traces and using the fact that $\alpha$ is
co-closed, we easily obtain the following relation:
\begin{equation}\label{fin}
\langle\mathrm{trace}_{g}(\tilde{B})(D^{2}A),A\rangle
=-8|A|^{2}\mathrm{trace}_{g}(\eta )
=-8|A|^{2}\left(\frac{1}{4(n+2)}\mathrm{Scal}^{g}-2|\alpha |^{2}\right) .
\end{equation}
On the other hand, from Proposition \ref{imp} and the very definition
of the Penrose operator,
$$
\tilde{B}(D^{2}A)=D\circ B(DA)=
4D^{2}A-6D(P^{D}A).
$$
In particular, if $A\in\Gamma (Q)$ is a solution of the Penrose
operator, then $P^{D}A=0$ and
\begin{equation}
\langle \mathrm{trace}_{g}(D^{2}A),A\rangle
=\frac{1}{4}\langle \mathrm{trace}_{g}(\tilde{B})(D^{2}A),A\rangle =
-2|A|^{2}\left(\frac{1}{4(n+2)}\mathrm{Scal}^{g}-2|\alpha |^{2}\right) .
\end{equation}
\end{proof}

We now restrict to the situation when $(M,Q,g)$ is a compact
quaternionic-K\"{a}hler manifold of negative or zero scalar curvature. 
An important class of compact 
quaternionic-K\"{a}hler manifolds 
of negative scalar curvature can be constructed in the
following way \cite{alek1}: take a non-compact symmetric
quaternionic-K\"{a}hler manifold M = G/K, which is dual to a Wolf
space. The non-compact simple Lie group G has a torsion free
co-compact discrete subgroup $\Gamma$. Then the double quotient
$M/\Gamma$ is a compact quaternionic-K\"{a}hler manifold of
negative scalar curvature.

\begin{thm}\label{doi} Let $D$ be a quaternionic connection on 
a compact quaternionic-K\"{a}hler manifold $(M, Q, g)$, 
such that $D=D^{g}+S^{\alpha}$, where $\alpha$ is co-closed.  
If $\mathrm{Scal}^{g}<0$ (respectively, 
$\mathrm{Scal}^{g}=0$ and $\alpha \neq 0$ at any point 
of a dense subset of $M$) then
$P^{D}$ has no non-trivial global solutions. 
\end{thm}

\begin{proof}
Choose an admissible basis $\{ J_{1},J_{2},J_{3}\}$ of $Q$ and
recall the formula
\begin{equation}\label{ra}
D_{X}A=D^{g}_{X}A+B(\alpha\otimes A)(X)=D^{g}_{X}A+\sum_{j=1}^{3}\alpha
\left( [J_{j},A](X)\right) J_{j},
\end{equation}
which relates $D$ and $D^{g}$ when
they act on the bundle $Q$ (above $A\in\Gamma (Q)$ and 
$X\in TM$). 
Using relation (\ref{ra}) it
is straightforward to check that
\begin {align*}
\langle (D^{2}A)_{X,X},A\rangle &=\langle (D^{g})^{2}(A)_{X,X},A\rangle
+\alpha\left([A,D^{g}_{X}A](X)\right)
-2\alpha (X)\langle D^{g}_{X}A,A\rangle\\
&+\sum_{j=1}^{3}[X\left( \alpha ([J_{j},A](X))\right)\langle J_{j},A\rangle
+\alpha\left( [J_{j},A](X)\right)\langle D_{X}J_{j},A\rangle]\\
&+2\sum_{j=1}^{3} \alpha (J_{j}X)\langle D^{g}_{J_{j}X}A,A\rangle .
\end{align*}
Choosing a $g$-orthonormal basis $\{ e_{1},\cdots ,e_{4n}\}$ of $TM$ and letting
$X:=e_{i}$ in the previous relation, we get
\begin{align*}
\langle\mathrm{trace}_{g}(D^{2}A), A\rangle &=
\langle\mathrm{trace}_{g}(D^{g})^{2}A, A\rangle
+4\langle D^{g}_{g^{*}\alpha}A,A\rangle +2\sum_{j=1}^{4n}
\alpha \left([A,D^{g}_{e_{j}}A]e_{j}\right)\\
&-\sum_{i,j=1}^{4n}\alpha\left( [J_{j},A]e_{i}\right)^{2}.
\end{align*}
On the other hand, again from relation
(\ref{ra}), we deduce that
\begin{align*}
\langle DA,DA\rangle &:=\sum_{i=1}^{4n}\langle D_{e_{i}}A,
D_{e_{i}}A\rangle =\sum_{i=1}^{4n}\left( \langle
D^{g}_{e_{i}}A,D^{g}_{e_{i}}A\rangle +2\alpha\left(
[D^{g}_{e_{i}}A,A]e_{i}\right)
\right)\\
&+\sum_{i,j=1}^{4n}\alpha \left( [J_{j},A]e_{i}\right)^{2}\\
&=\langle D^{g}A,D^{g}A\rangle +2\sum_{i=1}^{4n}\alpha\left(
[D^{g}_{e_{i}}A,A]e_{i}\right)+\sum_{i,j=1}^{4n}\alpha \left(
[J_{j},A]e_{i}\right)^{2}.
\end{align*}
Combining the above relations, we get
\begin{align*}
\langle\mathrm{trace}_{g}(D^{2}A),A\rangle =&
\langle\mathrm{trace}_{g}(D^{g})^{2}(A),A\rangle +4\langle
D^{g}_{g^{*}\alpha}A,A\rangle\\
&+\langle D^{g}A, D^{g}A\rangle -\langle DA,DA\rangle .
\end{align*}
Suppose now that $P^{D}A=0.$ Using Proposition \ref{propo}, this
relation becomes
\begin{align*}
&\langle DA,DA\rangle -\langle D^{g}A,D^{g}A\rangle-4\langle
D^{g}_{g^{*}\alpha}A,A\rangle -\langle
\mathrm{trace}_{g}(D^{g})^{2}A,A\rangle\\
&=2|A|^{2} \left(\frac{1}{4(n+2)}\mathrm{Scal}^{g}-2|\alpha |
^{2}\right) .
\end{align*}
Integrating over $M$ and using $\int_{M}\langle
D^{g}_{g^{*}\alpha}A,A\rangle\mathrm{vol}_{g} =0$, the $1$-form
$\alpha$ being co-closed, we get
\begin{equation}\label{our}
\int_{M}\langle DA,DA\rangle\mathrm{vol}_{g}
+4\int_{M}|A|^{2}|\alpha |^{2}\mathrm{vol}_{g}
-\frac{1}{2(n+2)}\mathrm{Scal}^{g}\int_{M}|A|^{2}\mathrm{vol}_{g}=0.
\end{equation}
Our claim readily follows from (\ref{our}).

\end{proof}

As an application of Theorem \ref{doi} we state:

\begin{cor}\label{co} Let $D$ be a closed quaternionic connection on a compact
quaternionic-K\"{a}hler manifold $(M,Q,g)$.
If $\mathrm{Scal}^{g}<0$ (respectively, $\mathrm{Scal}^{g}=0$ and 
$D$ is not exact) 
then $\Theta$ has no global non-trivial $\bar{\partial}^{D}$-holomorphic 
sections.\end{cor}

\begin{proof}
Let us consider an arbitrary $\bar{\partial}^{D}$-holomorphic
section $s$ of $\Theta .$ As in \cite{paul}, we can prove that
$s=\tilde{A}+{\mathcal J}\tilde{B}$, for two sections
$A,B\in\Gamma (Q).$ It can  be checked that
$\bar{s}=-\tilde{A}+{\mathcal J}\tilde{B}$ is also
$\bar{\partial}^{D}$-holomorphic, from where we deduce that both
$\tilde{A}$ and $\tilde{B}$ are $\bar{\partial}^{D}$-holomorphic.
Therefore, to prove our claim it is enough to show that there are
no global non-trivial $\bar{\partial}^{D}$-holomorphic
distinguished sections of $\Theta .$ The Levi-Civita connection
$D^{g}$ is exact, and hence $D=D^{g}+S^{\alpha}$, for some closed
$1$-form $\alpha\in\Omega^{1}(M).$ From Theorem \ref{main}, the
holomorphic structure $\bar{\partial}^{D}$ depends (up to
isomorphism) only on the cohomology class of $\alpha .$ Hence,
without loss of generality, we can take $\alpha$ to be harmonic.
Note that $\alpha \neq 0$ when $\mathrm{Scal}^{g}=0.$ 
Recall now that an harmonic 
non-trivial $1$-form on a connected Riemannian manifold cannot vanish 
identically
on an 
open subset of that manifold \cite{a}. 
We conclude from Proposition
\ref{lipsa} and Theorem \ref{doi}.
\end{proof}

{\bf Remark:} We now generalize Corollary \ref{co} to positive tensor
powers of the tangent vertical bundle $\Theta$
(this will also provide an alternative proof of Corollary
\ref{co}, which does not use the Penrose transform of Proposition 
\ref{lipsa}). More precisely,
in the setting of Corollary \ref{co} we prove that
any tensor power $\Theta^{s}$ (where
$s\in\mathbb{N}\setminus \{ 0\}$) admits no global non-trivial
$\bar{\partial}^{D}$-holomorphic sections. As above, we can suppose that 
$D=D^{g}+S^{\alpha}$ with $\alpha$ harmonic. 
On the twistor space $Z$ of $(M, Q)$ we define a
Riemannian metric $h_{g}^{D}$ by the following conditions: the
horisontal space $H^{D}$ determined by the connection $D$, acting
on $\pi :Z\to M$, is $h_{g}^{D}$-orthogonal to the twistor
lines; $h_{g}^{D}$, restricted to any horisontal subspace
$H^{D}_{J}$ (identified with $T_{\pi (J)}M$ by means of $\pi_{*}$), 
coincides with $g_{\pi (J)}$;
$h_{g}^{D}$ restricted to the twistor lines coincides with the
standard metric $\langle\cdot ,\cdot\rangle$ of the twistor lines.
Clearly, the pair $(h_{g}^{D}, {\mathcal J})$ is Hermitian. We are
interested in its torsion form $t_{g}^{D}\in\Omega^{1}(Z)$,
defined by
$$
t_{g}^{D}(\cdot ):=-\sum_{k
=1}^{4n+2}d\Omega_{g}^{D}(E_{k},{\mathcal J}E_{k},\cdot )
$$
where $\Omega_{g}^{D}:=h_{g}^{D}({\mathcal J}\cdot ,\cdot )$ is
the K\"{a}hler form and $\{ E_{1}, \cdots ,E_{4n+2}\}$ is a local
$h_{g}^{D}$-orthonormal frame of $Z$.

\begin{lem} The torsion form $t_{g}^{D}$ is equal to
$8\pi^{*}\alpha$ and is co-closed.
\end{lem}

\begin{proof} Like in Lemma 12 of \cite{paul}, we have the following expression
for $d\Omega_{g}^{D}$:
\begin{align*}
d\Omega_{g}^{D}(\tilde{X}, \tilde{Y},\tilde{V})&= (D_{X}g)(V, JY)
+(D_{Y}g)(X, JV)+(D_{V}g)(Y, JX)\\
d\Omega_{g}^{D}(a,\tilde{X},\tilde{Y})&= g(aX, Y)-\langle
[R^{D}_{X, Y}, J], Ja\rangle\\
d\Omega_{g}^{D}(a, b, c)&=0\\
d\Omega_{g}^{D}(a, b,\tilde{X})&=0,
\end{align*}
where $\tilde{X}, \tilde{Y},\tilde{V}\in H^{D}_{J}$ project to
$X,Y, V\in T_{p}Z$ (here $p:=\pi (J)$) and $a, b, c\in
T_{J}Z_{p}$. Since $D=D^{g}+S^{\alpha}$,
\begin{align*}
(D_{X}g)(Y, V)&=-2\alpha (X)g(Y, V)-\alpha (Y)g(X,V)-\alpha
(V)g(X,Y)\\
&+\sum_{i=1}^{3}[\alpha (J_{i}Y)g(J_{i}X, V)+\alpha
(J_{i}V)g(J_{i}X, Y)],
\end{align*}
where $\{ J_{1}, J_{2}, J_{3}\}$ is an admissible basis of $Q$; we
have similar expressions for $(D_{Y}g)(X, JV)$ and $(D_{V}g)(Y,
JZ).$ Using these observations, it is straightforward to show that
$t_{g}^{D}(\tilde{X})= 8\alpha (X).$ Moreover, if we chose the
admissible basis of $Q$ such that $J=J_{1}$, then any $a\in
T_{J}Z_{p}$ is of the form $a=\lambda_{2}J_{2}+\lambda_{3}J_{3}$
for $\lambda_{2}, \lambda_{3}\in \mathbb{R}$. Using relations
(\ref{alta1}) and (\ref{alta2}), we get:
$$
t_{g}^{D}(a)= \sum_{k =1}^{4n}\langle [R^{D}_{X_{k}, JX_{k}}, J],
Ja\rangle = \frac{1}{2n}\sum_{k=1}^{4n}\left(
\lambda_{2}\mathrm{tr}\left( J_{2}R^{\eta }_{X_{k},JX_{k}}\right)
+\lambda_{3}\mathrm{tr}\left( J_{3}R^{\eta }_{X_{k},JX_{k}}\right)
\right) ,
$$
where $R^{D}=W+R^{\eta}$ (see Section \ref{basicquat}) and $\{
X_{1},\cdots ,X_{4n} \}$ is a $g$-orthonormal basis of $T_{p}M.$
Using the definition of $R^{\eta}$ one can check that, for any
$i\in \{ 2,3\}$,
\begin{equation}\label{sums}
\sum_{k=1}^{4n}\mathrm{tr}\left( J_{i}R^{\eta}_{X_{k},
JX_{k}}\right) =8n \sum_{k=1}^{4n} \eta (JX_{k}, J_{i} X_{k}).
\end{equation}
Since $D$ is closed, $\eta$ is symmetric and the right hand side
of (\ref{sums}) is zero. Therefore, $t_{g}^{D}(a)=0$ and we can
conclude that $t_{g}^{D}=8\pi^{*}\alpha .$ Moreover, $t_{g}^{D}$
is co-closed since $\alpha$ is co-closed and $\pi :(Z,
h_{g}^{D})\to (M, g)$ is a Riemannian submersion with
totally geodesic fibers.
\end{proof}

We conclude the proof of our claim, by showing that
\begin{equation}\label{mic}
\int_{Z}\langle \gamma^{D},
\Omega_{g}^{D}\rangle_{h_{g}^{D}}\mathrm{vol}_{h_{g}^{D}}<0,
\end{equation}
where $\gamma^{D} := -\frac{1}{2\pi i}R^{\nabla}$ is the Chern
form of $\Theta$, endowed with its standard Hermitian structure
and Chern connection $\nabla$ induced by $D$, and $\langle\cdot
,\cdot\rangle_{h_{g}^{D}}$ is the inner product on $2$-forms (or
bivectors) on $Z$, defined, in our conventions,  by
$$
\langle V_{1}\land W_{1} ,V_{2}\land
W_{2}\rangle_{h_{g}^{D}}=h_{g}^{D}(V_{1}, V_{2}) h_{g}^{D}(W_{1},
W_{2})-h_{g}^{D}(V_{1}, W_{2})h_{g}^{D}(V_{2}, W_{1}).
$$
Since the pair $(h_{g}^{D},{\mathcal J})$ is standard  
(i.e. $t_{g}^{D}$ is co-closed) relation
(\ref{mic}) will insure that $\Theta^{s}$ (for any $s\in\mathbb{N}\setminus \{ 0\}$) 
has no
non-trivial global $\bar{\partial}^{D}$-holomorphic sections (see
\cite{paul1}, p. 504 and the argument used in \cite{paul}, p. 611). 
To prove (\ref{mic}), we notice that $\gamma^{D}$,
at a point $J\in Z_{p}$, has the following expression:
\begin{align*}
\gamma^{D}(\tilde{X}, \tilde{Y})&=\frac{1}{4n\pi
}\mathrm{tr}\left(
JR^{D}_{X, Y}\right)\\
\gamma^{D}(a, b)&=\frac{1}{2\pi}\Omega_{p}(a, b)\\
\gamma^{D}(\tilde{X}, a)&=0,
\end{align*}
where, we recall, $\Omega_{p}$ is the K\"{a}hler form of the
twistor line $Z_{p}$. The inner product $\langle \gamma^{D},
\Omega_{g}^{D}\rangle_{h_{g}^{D}}$ at $J$ is the
sum of the inner products of the vertical parts of $\gamma^{D}$
and $\Omega_{g}^{D}$, which is equal to $\frac{1}{2\pi }$, and of
their horisontal parts, which is equal to $\frac{1}{2} \sum_{k
=1}^{4n}\gamma^{D}(\widetilde{X_{k}}, \widetilde{JX_{k}}).$ We easily get that
\begin{align*}
\langle\gamma^{D},
\Omega_{g}^{D}\rangle_{h_{g}^{D}}=\frac{1}{2\pi}
+\frac{1}{4(n+2)\pi }\left(\mathrm{Scal}^{g}-8(n+2)|\alpha
|_{g}^{2}\right) .
\end{align*}
It follows that (\ref{mic}) holds if and only if
\begin{equation}\label{mic1}
\frac{1}{2\pi} +\frac{1}{4(n+2)\pi \mathrm{vol}_{g}(M)}
\int_{M}\left(\mathrm{Scal}^{g}-8(n+2)|\alpha|_{g}^{2}\right) \mathrm{vol}_{g}<0,
\end{equation}
where $\mathrm{vol}_{g}(M)=\int_{M}\mathrm{vol}_{g}$ denotes the
total volume of $M$ with respect to $g$. Since
$\mathrm{Scal}^{g}\leq 0$ and $\alpha \neq 0$ when $\mathrm{Scal}^{g}=0$,
we can make relation (\ref{mic1}) to hold,
by replacing, if necessary, $g$ with 
$tg$, where $t>0$ is sufficiently small. This concludes the proof
of our claim.\\

{\bf Acknowledgements:} I am the holder of a junior research
fellowship at Centro di Ricerca Matematica ``Ennio de Giorgi'',
Scuola Normale Superiore, Pisa. Hospitality and stimulating
working environment during my visit in May 2006 at Centre de
Mathematiques ``Laurent Schwartz'' of Ecole Polytechnique
(Palaiseau) is also acknowledged.

\end{document}